# Some numerical characteristic inequalities of compact Riemannian manifolds

Sergey Stepanov[1, 2 *], Irina Tsyganok[2]

[1] *Department of Scientific Information on Fundamental and Applied Mathematics, Russian Institute for Scientific and Technical Information of the Russian Academy of Sciences, Moscow, Russia*

[2] *Dept. of Mathematics and Data Analysis, Finance University, Moscow, Russia.*

**Abstract**

In this paper, we investigate certain numerical characteristics of an $n$-dimensional connected, compact Riemannian manifold $(M, g)$, such as the supremum and infimum of its Ricci curvature $Ric$ and sectional curvature $sec$, and explore their applications. We present two illustrative results. First, if $(M, g)$ is a compact, connected Riemannian manifold of even dimension $n = 2k \geq 4$ whose the supremum and infimum of its Ricci and sectional curvatures satisfy the strict inequality $2k\ sec_{\text{inf}} > Ric_{\text{sup}}$, then $M$ is diffeomorphic to either the Euclidean sphere $\mathbb{S}^{2k}$ of some radius $r > 0$ or the real projective space $\mathbb{RP}^{2k}$. Second, there exists no harmonic immersion of an $n$-dimensional compact, connected Riemannian manifold $(M, g)$ into the Euclidean sphere $\mathbb{S}^m$ of radius $r > 0$ if the infimum of its Ricci curvature satisfies the strict inequality $Ric_{\text{inf}} > n/2r^2$.

**Keywords** compact Riemannian manifold, supremum and infimum of the sectional and the Ricci curvatures. spherical theorem, harmonic map

## 1. Introduction

In Riemannian geometry, it is customary to consider only three classical numerical characteristics of a compact (boundaryless) Riemannian manifold $(M, g)$: its volume, diameter, and total scalar curvature. In this paper, we extend this perspective by focusing on additional numerical characteristics — namely, the supremum and infimum of the Ricci curvature $Ric$ and the sectional curvature $sec$ of the manifold. Our first main result is a new version of the Differentiable Sphere Theorem (see [1]) for compact, connected $n$-dimensional Riemannian manifolds. Using the convergence theorem for the Ricci flow and the inequality $n\ sec_{\text{inf}} > Ric_{\text{sup}}$ for the infimum and supremum of its sectional and Ricci curvatures we derive geometric and topological consequences that characterize such manifolds.

In addition, we prove new versions of classical vanishing theorems for harmonic

---

* Corresponding author.

E-mail address: s.e.stepanov@mail.ru (S.E. Stepanov)

maps. For instance, we show that there is no harmonic immersion of an $n$-dimensional compact, connected Riemannian manifold $(M, g)$ into a Euclidean sphere $\mathbb{S}^m$ of radius $r > 0$ if its Ricci curvature satisfies $Ric_{\inf} > n/2r^2$, where $Ric_{\inf}$ is the infimum of its Ricci curvature.

## 2. A Differentiable Sphere Theorem via the infimum and supremum of Ricci and sectional curvatures on compact manifolds

Let consider a compact (without boundary) Riemannian manifold $(M, g)$. We recall that the supremum and infimum of a continuous function $f: M \to \mathbb{R}$ defined on a compact manifold $M$ are the supremum and infimum of its range, and results about sets translate immediately to results about functions. Moreover, the extreme value theorem is well known, according to which a continuous function $f: M \to \mathbb{R}$ on a compact manifold $M$ attains its supremum and infimum, i.e., there exist $p, q \in M$ such that $f(p) = \text{Sup}\,\{f(x): x \in M\}$ and $f(q) = \text{Inf}\,\{f(x): x \in M\}$. Therefore, since the unit sphere in $T_xM$ at an arbitrary point $x \in M$ is a compact manifold, then there exist the unit vector $X_{\sup} \in T_xM$ and $X_{\inf} \in T_xM$ such that $Ric_{\sup}(x) \coloneqq Ric(X_{\sup})$ and $Ric_{\inf}(x) \coloneqq Ric(X_{\inf})$ for the real number $Ric_{\sup}(x) := \text{Sup}\,\{Ric(X): X \in T_xM\}$ and $Ric_{\inf}(x) := \text{Inf}\,\{Ric(X): X \in T_xM\}$, respectively, where $Ric(X)$ is the Ricci curvature in the direction to any unit vector $X \in T_xM$ at $x \in M$. Next, since $(M, g)$ is compact, then there exist the infinium and supremum of its Ricci $Ric$ curvature. Let us define and denote these numbers as follows: $Ric_{\inf} := \text{Inf}\,\{Ric_{\inf}(x): x \in M\}$ and $Ric_{\sup} := \text{Sup}\,\{Ric_{\sup}(x): x \in M\}$. Let us continue similar reasoning to determine the infimum of the section curvature on a compact manifold $(M, g)$. We first denote by $sec(\pi_x)$ the sectional curvature of $(M, g)$ in the direction to the 2-plane $\pi_x \subset T_xM$ at $x \in M$. At the same time, since the unit sphere in $T_xM$ at an arbitrary point $x \in M$ is a compact manifold, then there exist the 2-plane $\pi_{inf} = span\,\{X_{\inf}, Y_{\inf}\} \subset T_xM$ for unit orthogonal vectors $X_{\inf}, Y_{\inf} \in T_xM$ at $x \in M$ such that $sec_{\inf}(x) := sec_{\inf}(\pi_{inf})$, where $sec_{\inf}(x) := \text{Inf}\,\{sec(\pi_x): \pi_x \subset T_xM\}$. Next, since $(M, g)$ is compact, then there exists the infinium of its sectional curvature. Let us designate and define this number as $sec_{\inf} := \text{Inf}\,\{sec_{\inf}(x): x \in M\}$. Additionally, note that $sec_{\inf}$ is another numerical invariant (along with $Ric_{\inf}$ and $Ric_{\sup}$) of a compact manifold $(M, g)$, determined by its metric $g$.

On the other hand, we recall here that the Differentiable Sphere Theorem asserts (see [1, p. 4]) that any compact, connected Riemannian manifold $(M, g)$ which is strictly $1/4$-pinched in the pointwise sense admits another Riemannian metric which has constant sectional curvature. In particular, this implies that $M$ is diffeomorphic to a spherical space form $\mathbb{S}^n/\Gamma$, where $\mathbb{S}^n$ is the sphere of some radius $r > 0$ in the Euclidean space $\mathbb{R}^{n+1}$ with the metric $g_{can}$ induced from the Euclidean metric on $\mathbb{R}^{n+1}$ and $\Gamma$ is a discrete group of isometries acting properly discontinuously on $\mathbb{S}^n$. In turn, we formulate a new version of the Differentiable Sphere Theorem.

**Theorem 1**. *Let $(M, g)$ be an $n$-dimensional $(n \geq 3)$ connected compact Riemannian manifold $(M, g)$ the inequality $n\, sec_{\inf} > Ric_{\sup}$ holds, where $sec_{\inf}$ and $Ric_{\sup}$ are the infimum and supremum of its sectional and Ricci curvatures, respectively. Then $M$ is diffeomorphic to a spherical space form $\mathbb{S}^n/\Gamma$. In particular, if $M$ is simply connected then it is diffeomorphic to $\mathbb{S}^n$.*

*Proof.* First, we recall that the Riemann curvature tensor $Rm$ of $(M, g)$ induces a self-adjoint with respect to the pointwise inner product on the space of trace-less symmetric two-tensors $S_0^2(T_xM)$. This operator, denoted $\mathring{R}: S_0^2(T_xM) \to S_0^2(T_xM)$, is called the *curvature operator of the second kind* at each point $x \in M$. We would like to mention that the action of the Riemann curvature tensor on trace-less symmetric two-tensors has a long and interesting history. For example, Nishikawa in [5], conjectured that compact manifolds with positive curvature operators of the second kind are diffeomorphic to spherical space forms $\mathbb{S}^n/\Gamma$.

In turn, Cao-Gursky-Tran in [6] proved Nishikawa's conjecture. Their key observation is that the two-positive curvature operator of the second kind implies a strictly PIC1 condition introduced by Brendle [7]. The positive case of Nishikawa's conjecture follows immediately from Brendle's result in [7] asserting that the normalized Ricci flow evolves an initial metric satisfying strictly PIC1 into a limit metric with constant positive sectional curvature for $n \geq 4$.

Second, from Theorem 1.1 of our paper [8] we know that in dimensional $n \geq 3$ a compact Riemannian manifold $(M, g)$, whose sectional and Ricci curvature satisfy the inequality $n\, sec_{\inf}(x) > Ric_{\sup}(x)$ at each point $x \in M$ has positive curvature operator of the second kind. Therefore, an $n$-dimensional, $n \geq 4$, compact manifold $M$ diffeomorphic to spherical space forms $\mathbb{S}^n/\Gamma$, if $n\, sec_{\inf}(x) > Ric_{\sup}(x)$ at each point $x \in M$.

Third, let the inequality $n\, sec_{\inf} > Ric_{\sup}$ holds. At the same time, it is clear that $n\, sec(\pi_x) \geq n\, sec(\pi_x) > Ric_{\sup} \geq Ric_{\sup}(x)$ for an arbitrary unit vector $X \in T_xM$ at each point $x \in M$, then the inequality $n\, sec_{\inf}(x) > Ric_{\sup}(x)$ holds. Therefore, a metric $g$ of a compact Riemannian manifold $(M, g)$ has positive sectional curvature if the inequality $n\, sec_{\inf} > Ric_{\sup}$ holds. In particular, for the Euclidean sphere $(\mathbb{S}^n, g_{can})$ of some radius $r > 0$ the equality $n\, sec_{\inf} > Ric_{\sup}$ is satisfied automatically.

And finally, fourthly, we can write the following inequalities:

$$s(x) \geq n(n-1)\, sec_{\inf}(x) > (n-1)\, Ric_{\sup}(x),$$

for the scalar curvature $s(x) := \sum_{i \neq j} sec(e_i, e_j)$ of $(M, g)$, where $sec(e_i, e_j)$ denotes the sectional curvature of $(M, g)$ at an arbitrary point $x \in M$ with respect to the 2-plane $\pi_x = span\,\{e_i, e_j\} \subset T_xM$ and $\{e_1, \dots, e_n\}$ is an orthonormal basis for $T_xM$. Thus from the inequality $n\, sec_{\inf} > Ric_{\sup}$ we obtain the inequality $s(x) > (n-1)\, Ric_{\sup}(x)$. In particular, for the case when $n = 3$ this inequality can be rewritten as $s(x) > 2\, Ric_{\sup}(x)$. At the same time, from Corollary 8.2 of the paper [4, p. 277] we know that a metric $g$ of three-dimensional Riemannian manifold $(M, g)$ has positive sectional curvature if $s(x)\, g(x) > 2Ric(x)$ at each point $x \in M$. Moreover, from Main Theorem 1.1 of the paper [4, p. 255] we know the following: A three-dimensional compact Riemannian manifold $(M, g)$ with positive Ricci curvature is diffeomorphic to a spherical space form $\mathbb{S}^3/\Gamma$. In particular, if $M$ is simply connected then it is diffeomorphic to $\mathbb{S}^3$. As a result, we conclude that Theorem 1 is true.

We will prove two consequences of Theorem 1 in conclusion of this section. The first corollary is the following.

**Corollary 1.** *Let* $(M, g)$ *be an* $2k$*-dimensional* $(k \geq 2)$ *connected, compact Riemannian manifold whose Ricci and sectional curvatures satisfy the strict inequality* $2k\, sec_{\inf} > Ric_{\sup}$*, where* $sec_{\inf}$ *and* $Ric_{\sup}$ *are their infimum and supremum, respectively. Then M is diffeomorphic to a spherical space form* $\mathbb{S}^{2k}/\Gamma$*, which is the sphere* $\mathbb{S}^{2k}$ *or the real projective space* $\mathbb{RP}^{2k}$.

*Proof.* From Theorem 2 we known ерy ащддщццштп: If $(M, g)$ is an $n$-dimensional $(n \geq 3)$ connected compact Riemannian manifold $(M, g)$ *and* the inequality $n\, sec_{inf} > Ric_{sup}$ holds, where $sec_{inf}$ and $Ric_{sup}$ are the infimum and supremum of its sectional and Ricci curvatures, respectively, then $M$ is diffeomorphic to a compact spherical space form $\mathbb{S}^n/\Gamma$. At the same time, the simplest examples of the spherical space form are the sphere $\mathbb{S}^n$ and the real projective space $\mathbb{RP}^n$. Furthermore, $n$ is even, these are the only examples (see [1, p. 3]). Therefore, Corollary 1 is a consequence of Theorem 2 proved above.

**Corollary 2.** *Let* $(M, g)$ *be an* $n$*-dimensional* $(n \geq 3)$ *connected compact Riemannian manifold* $(M, g)$ *whose Ricci tensor and sectional curvature satisfy the strict inequality* $n\, sec_{\inf} > Ric_{\sup}$*, where* $sec_{\inf}$ *and* $Ric_{\sup}$ *are infimum and supremum of its sectional and Ricci curvatures, respectively. Then there is no non-zero traceless Codazzi* $p$*-tensor* $(p \geq 2)$.

*Proof.* In the article [9], the authors considered the concept of Codazzi $p$-tensors $(p \geq 2)$ which extends the well-known concept for $p = 2$ (see [10, pp. 436-440]). Let us recall that a Codazzi $p$-tensor $(p \geq 2)$ or, in other words, a higher order Codazzi tensor is a $C^\infty$-section $\varphi$ of the vector bundle $S^pM$ of symmetric $p$-forms $(p \geq 2)$ on *M* satisfying the following condition: $\nabla\varphi \in C^\infty(S^{p+1}M)$. We proved in [11] that every traceless Codazzi $p$-tensor $\varphi \in C^\infty\left(S_0^pM\right)$ on a compact Riemannian manifold $(M, g)$ with nonnegative curvature operator of the second kind is invariant under parallel translations, i.e., $\nabla\varphi = 0$. Moreover, if $\mathring{R} > 0$ at some point $x \in M$, then there is no non-zero traceless Codazzi $p$-tensor ($p \geq 2$). Therefore, the second corollary holds.

## 3. Inequalities involving the infimum and supremum of Ricci and sectional curvatures and their applications to harmonic map theory

A systematic study of harmonic maps was initiated in 1964 by Eells and Sampson [2]. Detailed presentations of the results can be found in [13] – [16], [20] and many other publications including monographs (see, for example, [3]). For definitions, notations, and results, we will refer to these works.

Let $(M, g)$ be an $n$-dimensional ($n \geq 2$) connected Riemannian manifold with the Levi-Civita connections $\nabla$ and $(\bar{M}, \bar{g})$ be an $m$-dimensional ($m \geq 2$) connected Riemannian manifold $(\bar{M}, \bar{g})$ with the Levi-Civita connections $\bar{\nabla}$, and $f: (M, g) \to (\bar{M}, \bar{g})$ be a smooth mapping. Denote by $f_*$ the differential of the mapping $f$ and by $f^*$ we denote the dual mapping; $f^*T\bar{M}$ is a bundle over $M$ with fiber $T_{f(x)}\bar{M}$ over a point $x \in M$. Therefore, the differential $f_*$ is a section of the tensor bundle $f^*T\bar{M} \otimes T^*M$. The metric $\bar{g}$ and the connection $\bar{\nabla}$ on $(\bar{M}, \bar{g})$ induce the metric tensor $\bar{g}^* := f^*\bar{g}$ which is called in [14, p. 8] as the *first fundamental form* of $f$ and the connection $\bar{\nabla}^*$ on the bundle $f^*T\bar{M}$. Denote by $D$ the connection in the bundle $f^*T\bar{M} \otimes T^*M$ induced by the connections $\bar{\nabla}^*$ and $\nabla$ (see [14, p. 6]). Then for any vector fields $X, Y$ on $M$ we have $(D\, f_*)(X, Y) = (D_X f_*)Y = \bar{\nabla}^*_X(f_*Y) - f_*\nabla_X Y$, where $f_*Y$ is differentiated as a vector field along the mapping $f$ (see, for example, [14]). At the same time, the symmetric bilinear form $Df_*: TM \times TM \to f^*T\bar{M}$ is called the *second fundamental form* of the mapping $f$ (see [1, p. 2]; [14, p. 8]; [22, p. 389]). A map with vanishing second fundamental form is said to be *totally geodesic* (see [14, p. 9]; [22, p. 389]). Such a map is characterized by the property that it carries geodesies to geodesies linearly.

The *Dirichlet energy functional* of a smooth map $f: (M, g) \to (\bar{M}, \bar{g})$ over any compact domain $U \subset M$ is defined to be functional

$$E(f, U) = \int_U e(f)\, dv_g,$$

where $dv_g$ denotes the volume element of $(M, g)$ and $e(f) := (trace_g\,(f^*\bar{g}))$ is a non-negative scalar function (see [14, p. 10]). The term $e(f)$ is known (see [14, p. 10]) as the *energy density* of $f$, and it provides a measure of how much the map $f$ distorts or stretches the metric $\bar{g}$ of the target manifold $\bar{M}$ at each point in $M$. For an isometric immersion $f: (M, g) \to (\bar{M}, \bar{g})$, the identity $e(f) = n/2$ holds (see [14, p. 10]). Note also that $e(f) = 0$ if and only if $f$ maps $M$ to a point of $\bar{M}$, i.e., $f$ is a constant map.

*Euler-Lagrangian* associated with the Dirichlet energy functional is the *tension field* $\tau(f)$ that given by $\tau(f) = trace_g(Df_*)$ (see [14, p. 9]). Therefore, a smooth map $f: (M, g) \to (\bar{M}, \bar{g})$ is a *harmonic map* if and only if its tension field $\tau(f)$ vanished identically (see the definition, properties and examples of harmonic map in [1]; [3] [14] and [22]).

On the other hand, along with harmonic mappings, there are biharmonic mappings as a special case of $k$-polyharmonic mappings, first proposed by J. Eells and L. Lehmer in [23]. A biharmonic map $f: (M, g) \to (\bar{M}, \bar{g})$ between Riemannian manifolds is a critical point of the *Dirichlet bienergy functional* (see [25])

$$E_2(f, U) = \int_U \|\tau(f)\|^2\, dv_g$$

for any compact domain $U \subset M$. *Euler-Lagrangian* associated with the Dirichlet bienergy functional is the *Jacobi operator* $J(f)$ that given by $J\big(\tau(f)\big) = \Delta_f \tau(f) - trace\ \bar{R}\big(f_*,\ \tau(f)\big)f_*$, where $\Delta_f$ is the Laplacian on section $f^*T\bar{M}$ (see [23, p. 28]). Therefore, a smooth map $f: (M, g) \to (\bar{M}, \bar{g})$ is said to be biharmonic if $J\big(\tau(f)\big) = 0$ (see [25, p. 4]). At the same time, if a Riemannian manifold $(M, g)$ is compact and $(\bar{M}, \bar{g})$ is a Riemannian manifold with non-positive sectional curvature, then $f$ is biharmonic if and only if it is harmonic (see [24]).

For the case when $(M, g)$ is a compact manifold we deduce in [21] the integral formula

$$\int_M (Q(f) + \|Df_*\|^2 - \|\tau(f)\|^2) dv_g = 0, \qquad (1)$$

where $\|Df_*\|^2$ is the square of the norm of the second fundamental form of $f$ with respect to the metric on $T^*M\otimes T^*M\otimes f^*T\bar{M}$ induced by the metrics $g$ and $\bar{g}^*$. In turn, the scalar function $Q(f)$ has the form (see, for example, [13, p. 3]; [14, p. 12])

$$Q(f) = trace_{\bar{g}}(f_*Ric) - \sum_{i,j=1}^{n}(f^*\bar{R})(e_i, e_j, e_j, e_i) \tag{2}$$

where $Ric$ is the Ricci tensor of $(M, g)$, $\bar{R}$ is the Riemannian curvature tensor of $(\bar{M}, \bar{g})$ and $\{e_1, \ldots, e_n\}$ is the local orthonormal frame of $TM$. Moreover, the inequality

$$trace_{\bar{g}}(f_*Ric) \geq \frac{1}{n} Ric_{inf} \cdot e(f)$$

holds for the infimum $Ric_{inf}$ of Ricci curvature $Ric$ of a compact $(M, g)$.

In addition, if $(M, g)$ is compact, then we can define the *Dirichlet energy* and *bienergy* of the smooth mapping $f: (M, g) \to (\bar{M}, \bar{g})$ as its non-negative numeric characteristics $E(f) = \int_M e(f)\, dv_g$ (see [22, p. 388)] and $E_2(f) = \int_M \|\tau(f)\|^2\, dv_g$ (see [24]) respectively. At the same time, we can rewrite the last formula as $E_2(f) = \int_M (Q(f) + \|Df_*\|^2) dv_g$.

In particular, if $f: (M, g) \to (\bar{M}, \bar{g})$ is an isometric immersion, then $E(f) = n/2\ Vol(M) > 0$. Moreover, $E(f) = 0$ if and only if $f$ is a constant map since $E(f) = \int_M e(f)\, dv_g$ for $e(f) \geq 0$. On the other hand, $E_2(f) = 0$ if and only if $f$ is a harmonic map since $E_2(f) = \int_M \|\tau(f)\|^2\, dv_g$ for $\|\tau(f)\|^2 \geq 0$.

From the above we deduce the following:

$$0 = \int_M (Q(f) + \|Df_*\|^2 - \|\tau(f)\|^2) dv_g \geq \int_M Q(f)\, dv_g - E_2(f) \geq$$
$$\geq \int_M \left(-\sum_{i,j=1}^{n}(f^*\bar{R})(e_i, e_j, e_j, e_i)\right) dv_g + \frac{1}{n} Ric_{inf} \cdot E(f) - E_2(f). \tag{3}$$

If $(\bar{M}, \bar{g})$ is a well-known *Hadamard manifold* (see, for example, [17]), namely, a complete and noncompact, simply connected Riemannian manifold of non-positive sectional curvature and the equality $Ric_{inf} \cdot E(f) > n\, E_2(f)$ holds, then from (3) we deduce that $f$ must be a constant map. As a result, we have the following theorem.

**Theorem 2**. *Let* $f: (M, g) \to (\bar{M}, \bar{g})$ *be a smooth map between an n-dimensional* $(n \geq 2)$ *connected, compact Riemannian manifold* $(M, g)$ *and an* $m$*-dimensional* $(m \geq 2)$ *Riemannian manifold* $(\bar{M}, \bar{g})$ *of non-positive sectional curvature* $\overline{sec} \leq 0$ (*or, in particular,* $(\bar{M}, \bar{g})$ *be a Hadamard manifold*). *If the Dirichlet energy* $E(f)$ *and bienergy* $E_2(f)$ *of* $f$ *satisfy the inequality*

$$E_2(f) < \frac{1}{n} Ric_{inf} \cdot E(f),$$

*where* $Ric_{inf}$ *is the infimum of Ricci curvature of* $(M, g)$, *then* $f$ *is a constant map.*

In particular, if $f$ is a harmonic map, then from (1) we obtain

$$\int_M (Q(f) + \|Df_*\|^2) dv_g = 0. \quad (4)$$

Therefore, suppose that the Ricci curvature of $(M, g)$ is non-negative and the sectional curvature of $(\bar{M}, \bar{g})$ is non-positive, then from (2) we obtain $Q(f)(x) \geq 0$. Therefore, we can state the following: if $f: (M, g) \to (\bar{M}, \bar{g})$ is any harmonic map between a compact Riemannian manifold $(M, g)$ with non-negative Ricci tensor and a Riemannian manifold $(\bar{M}, \bar{g})$ with non-positive sectional curvature, then $f$ is *totally geodesic* and has constant the energy density $e(f)$. Furthermore, if there is at least one point of $(M, g)$ at which its Ricci curvature is positive, then every harmonic map $f: (M, g) \to (\bar{M}, \bar{g})$ is constant.

When $(\bar{M}, \bar{g})$ is a well-known Hadamard manifold, then the following corollary follows from the above theorem, the theorem from [24] and our Theorem 2.

**Corollary 3**. *Let* $f: (M, g) \to (\bar{M}, \bar{g})$ *be a biharmonic map between an n-dimensional* $(n \geq 2)$ *connected, compact Riemannian manifold* $(M, g)$ *with the positive infimum* $Ric_{inf}$ *and an* m*-dimensional* $(m \geq 2)$ *Riemannian manifold* $(\bar{M}, \bar{g})$ *of non-positive sectional curvature* $\overline{sec} \leq 0$ (*or, in particular,* $(\bar{M}, \bar{g})$ *be a Hadamard manifold*). *Then* $f$ *is a constant map.*

Now, we formulate and prove an analogue of the above Eells-Sampson's theorem on harmonic mappings between Riemannian manifolds.

**Theorem 3** *Let* $f: (M, g) \to (\bar{M}, \bar{g})$ *be a harmonic mapping between connected, compact Riemannian manifolds such that* $\dim M = n \geq 2$ *and* $\dim M = m \geq 2$. *Suppose that the second manifold* $(\bar{M}, \bar{g})$ *satisfies the curvature inequality* $m\, \overline{sec}_{\inf} > \overline{Ric}_{\sup}$, *where* $\overline{sec}_{\inf}$ *and* $\overline{Ric}_{\sup}$ *denote the infimum the sectional curvature and the supremum of the Ricci curvature of* $(\bar{M}, \bar{g})$, *respectively.*

*Furthermore, assume that the energy density* $e(f)$ *of* $f$ *satisfies the inequality* $e(f)_{\sup} < \mathrm{Ric}_{\inf} / \overline{sec}_{\inf}$*, where* $\mathrm{Ric}_{\inf}$ *is the infimum of the Ricci curvature of* $(M, g)$. *Then* $f$ *is a constant map.*

*Proof.* We recall that if $f$ is a harmonic map, then the Bochner-type formula

$$\int_M \left(\|Ddf\|^2 + Q(f)\right) dv_g = 0 \tag{6}$$

holds. In this case, we obtain from (6) the following integral inequality (see [1, p. 124])

$$\int_M Q(f)\, dv_g \leq 0. \tag{7}$$

On the other hand, we proved the following equality (see [16])

$$Q(f) = \sum_{\alpha<\beta} \overline{sec}\left(\bar{e}_\alpha, \bar{e}_\beta\right)\left(\lambda_\alpha - \lambda_\beta\right)^2 + g(Ric - f^*\overline{Ric}, f^*\bar{g})$$

where $\overline{sec}\left(\bar{e}_\alpha, \bar{e}_\beta\right)$ is the sectional curvature of $(\bar{M}, \bar{g})$ in direction to $\bar{\pi}_{f(x)} = \mathrm{span}\left\{\bar{e}_\alpha, \bar{e}_\beta\right\} \subset T_{f(x)}\bar{M}$ at an arbitrary point $f(x) \in \bar{M}$ for $\lambda_\alpha \geq 0$ and $\alpha, \beta = 1, \ldots, m$. At the same time, we proved in [18] that $\sum_{\alpha=1}^{m}(\lambda_\alpha)^2 = \|f^*\bar{g}\|^2$ and $\sum_{\alpha=1}^{m} \lambda_\alpha = e(f)$ at an arbitrary point $f(x) \in \bar{M}$ . Since

$$\sum_{\alpha<\beta}\left(\lambda_\alpha - \lambda_\beta\right)^2 = m \sum_\alpha (\lambda_\alpha)^2 - \left(\sum_\alpha \lambda_\alpha\right)^2 = m\,\|f^*\bar{g}\|^2 - e(f)^2,$$

the following inequalities hold:

$$Q(f) \geq \overline{sec}_{\inf}(x)\left(m \sum_\alpha (\lambda_\alpha)^2 - \left(\sum_\alpha \lambda_\alpha\right)^2\right) + g(Ric - f^*\overline{Ric}, f^*\bar{g}) \geq$$

$$\geq \overline{sec}_{\inf}(x)\left(m\,\|f^*\bar{g}\|^2 - e(f)^2\right) + Ric_{\inf}(x)\, e(f) - \overline{Ric}_{\sup}(x) \cdot \|f^*\bar{g}\|^2 =$$

$$\geq \|f^*\bar{g}\|^2\left(m\,\overline{sec}_{\inf} - \overline{Ric}_{\sup}\right) + e(f)\left(Ric_{\inf} - e(f) \cdot \overline{sec}_{\inf}\right). \tag{8}$$

Moreover, if the inequality $\overline{sec}_{\inf} > 1/m\ \overline{Ric}_{\sup}$ holds, then $(\bar{M}, \bar{g})$ has positive sectional curvature, in contrast to Eells-Samson's theorem, and $\bar{M}$ is diffeomorphic to a spherical space form $\mathbb{S}^m/\Gamma$ (see our Theorem 2). We denote by $e(f)_{\sup} := \mathrm{Sup}\,\{e(x) : x \in M\}$. In this case, if the inequality $\overline{sec}_{\inf} \cdot e(f)_{\sup} < Ric_{\inf}$ holds, then $(M, g)$ has positive Ricci curvature. Furthermore, from the above we obtain the following integral inequality:

$$\int_M Q(f)\, dv_g = \left(m\,\overline{sec}_{\inf} - \overline{Ric}_{\sup}\right) \int_M \|f^*\bar{g}\|^2\, dv_g +$$

$$+ \int_M e(f)\left(Ric_{\inf} - e(f) \cdot \overline{sec}_{\inf}\right) dv_g \geq 0,$$

which contrasts with inequality (7). This forces $f$ to be a constant map.

In conclusion of our article, we will formulate two consequences of Theorem 3. The first of these consequences is obvious.

**Corollary 4.** *Let $f:(M,g)\to(\mathbb{S}^m,g_{can})$ be a harmonic map between a connected compact $n$-dimension Riemannian manifold $(M,g)$ of positive Ricci curvature and the $m$-dimension Euclidean sphere $(\mathbb{S}^m,g_{can})$ of radius r. If the energy density $e(f)$ of $f$ satisfies the inequality $e(f)_{\sup} < Ric_{\inf}$, then $f$ is a constant map.*

Consider a biharmonic immersion $f:(M,g)\to(\bar{M},\bar{g})$ of a compact Riemannian manifold $(M,g)$ into a compact Riemannian manifold $(\bar{M},\bar{g})$ (see [1, pp. 124-125[). In this case, the equalities $e\,(f) = n/2$ and $Q(f) + \|Ddf\|^2_{\bar{g}} = 0$ hold. In this case (8) can be rewritten as

$$Q(f) \geq \|f^*\bar{g}\|^2\left(m\,\overline{sec}_{\inf} - \overline{Ric}_{\sup}\right) + \frac{n}{2}\left(Ric_{\inf} - \frac{n}{2}\overline{sec}_{\inf}\right).$$

Therefore, if the double inequality $1/m\ \overline{Ric}_{\sup} < \overline{sec}_{\inf} < 2/n\ Ric_{\inf}$ holds, then $f:(M,g)\to(\bar{M},\bar{g})$ can not be a biharmonic immersion of $(M,g)$ into $(\bar{M},\bar{g})$ since $Q(f) > 0$. As a result, we can formulate the second corollary.

**Corollary 5.** *Let $(M,g)$ and $(\bar{M},\bar{g})$ be a connected compact Riemannian manifolds such that $1/m\ \overline{Ric}_{\sup} < \overline{sec}_{inf} < 2/n\ Ric_{\inf}$, then there is no harmonic immersion of $(M,g)$ into $(\bar{M},\bar{g})$.*

The celebrated theorem of Myers (see [19]) guarantees the compactness of complete Riemannian manifolds under some positive lower bounds on the Ricci curvature: Let $(M,g)$ be an $n$-dimensional complete Riemannian manifold. Suppose that there exists some positive constant $\lambda > 0$ such that the Ricci curvature satisfies the inequality $Ric \geq \lambda g$. Then $(M,g)$ must be compact with finite fundamental group. Let us reformulate this theorem as follows: If there exists some positive constant $\lambda > 0$ such that the Ricci curvature satisfies the inequality $Ric_{\inf} \geq \lambda$, then complete manifold $(M,g)$ must be compact with finite fundamental group. Using this proposition and corollaries 4 and 5 we can formulate our last statement.

**Corollary 6.** *There is no biharmonic immersion of a connected complete Riemannian manifold* $(M, g)$ *into the Euclidean sphere* $(\mathbb{S}^m, g_{can})$ *of radius* $r$ *if there exists* $Ric_{\inf}$ *such that* $Ric_{\inf} > n/2r^2$.

**Conflict of interest.** The authors declare that they have no conflict of interest.

**Funding Declaration.** The authors declare that no funds, grants, or other support were received during the preparation of this manuscript.